\documentclass[12 pt]{article}
\usepackage{graphicx}
\pdfoutput=1
\usepackage{amsmath}
\usepackage[ruled,vlined]{algorithm2e}

\usepackage{dcolumn}
\usepackage{amsmath} 
\usepackage{amssymb} 
\usepackage{enumerate}
\usepackage{amsthm}
\usepackage{float} 
\usepackage{capt-of} 
\usepackage{sidecap} 
\sidecaptionvpos{figure}{c} 
\usepackage{caption} 
\usepackage{commath} 
\usepackage{cancel} 
\usepackage{anysize} 
\usepackage{appendix} 
\usepackage{tocbibind}
\usepackage{physics} 
\usepackage{multirow}
\usepackage{glossaries}

\begin{document}
	
\title{Positive invariant regions for a modified Van Der Pol equation modeling heart action}
\author{
A. Acosta, R. Gallo\\
\small School of Mathematics, Yachay Tech, Yachay City of Knowledge, Ecuador.\\

P. Garc\'{\i}a\\
\small		Laboratorio de Sistemas Complejos,
\small		Departamento de F\'{\i}sica Aplicada,\\	
\small		Facultad de Ingenier\'{\i}a, 	\small		Universidad Central de Venezuela.\\	
\small		Red Iberoamericana de Investigadores en
\small		Matemáticas Aplicadas a  Datos. AUIP.\\

D. Peluffo-Ord\'o\~nez\\
\small Modeling, Simulation and Data Analysis (MSDA) Research Program,\\
\small  Mohamed VI Polytechnic University, Morocco.\\
}

	\date{}
	
	\maketitle
	\thispagestyle{empty}
\begin{abstract}
	The dynamical characterization of the heart rate is definitely a problem of vital importance. The selection, construction and adjustment of models that reproduce the dynamic behavior of the cardiac muscle, brings us closer to the solution of the usual classification problem in medicine, i.e. decide whether or not a patient belongs to the healthy class of patients.
	
	An usual model for this dynamic is often given by a  modified Van der Pol model (vdPM), whose parameters are currently optimized with machine learning techniques or evolutionary algorithms. In any of these cases, the proper selection of the initial condition of the parameters drastically reduces the computational cost of the optimization method used.
	
	In this work, a strategy to estimate a positively invariant region that contains the periodic orbit, associated with a set of particular parameters of the vdPM considered,  is presented. An algorithm is proposed to build the positively invariant region, the numerical results confirm that the vdPM solutions starting within the positively invariant region converge to the periodic orbit. This kind of result allows to reduce the cost of searching for the optimal parameters that represent a real electrocardiographic signal, which in turn, let us reduce the dimensionality of the system allowing the design of more efficient classifiers for these signals.
\end{abstract}

\section{Introduction}
The extensive work done by van der Pol in the decade of the 20s led to what is known today as the classical van der Pol oscillator and its usual presentation is given by the second order ordinary differential equation
\begin{eqnarray}
\ddot{x} + \alpha ( x^{2}  - 1 )  \dot{x} + x& = & 0  \label{original} 
\end{eqnarray}
where $  \alpha  $ is a scalar parameter indicating the strenght of the nonlinear damping. On the other hand, the first model of the heart based on equation (\ref{original}) was proposed in 1928 by van der Pol and van der Mark \cite{VDP2} and since then, inspired in equation (\ref{original}), much attention has been paying to the study  and formulation of heart action models; see for instance \cite{Acosta}, \cite{Zduniak}, \cite{Fitz}, \cite{West}, \cite{Bernardo}, \cite{VDPM}, \cite{Ryz}, \cite{RicardoViana}, \cite{Zebroski2}. In references \cite{Acosta}, \cite{Zduniak}, \cite{VDPM}, \cite{Zebroski2} the authors consider modifications of (\ref{original}) by adding two fixed points: a saddle and a node. This modified van der Pol model was capable to reproduce, in a better way that previous models based in equation (\ref{original}), the time series of the action potential generated by a natural pacemaker of the heart. Specifically in \cite{Acosta}, \cite{VDPM}, \cite{Zebroski2} the authors  consider the equation
\begin{eqnarray}
\ddot{x} + \alpha ( x^{2}  - {\nu}^{2} )  \dot{x} + \frac{x (x + d ) ( x + e )}{ed}& = & 0 ~ , \label{lamodificada} 
\end{eqnarray}
where $  \alpha , {\nu} , d , e  $  are positive control parameters.  
Because we are seeing the equation (\ref{lamodificada}) as a representation of a model that studies the actions of the heart, it is clear that the identification of regions where periodic orbits appear is of fundamental importance. In that sense, as far as we know, not much has been done and a case we have found in the literature, \cite{Zduniak}, corresponds to a specific configuration of the parameters for which the authors exhibit a positively invariant region, which was found using numerical interpolation, which captures periodic orbits of interest for the model. Thus, in general, researchers mostly focus on finding configurations of the parameters that produce periodic orbits and the search is carried out through the implementation of genetic algorithms (GA) \cite{Acosta}, optimization by particle swarm particles\cite{Lopes}, neural networks (NN)\cite{Khan}, mixed NN and GA\cite{Raja} or some other numerical strategy, see for example \cite{Suraj}.

\noindent
In this work we focus, mainly, on the identification of regions that capture periodic orbits for the equation  (\ref{lamodificada}). 

We construct a positively invariant region bounded by a closed curve that will encircles  the origin. To accomplish this goal, we convert the equation (\ref{lamodificada}) to an equivalent planar system of first-order equations. Thus, the  closed curve will be obtained by putting together pieces of orbits coming from planar systems that result when neglecting terms on the equivalent system, pieces of nullclines of the the equivalent system and some line segments.

Our work is a particularity of a more general problem which consists of constructing closed curves, depending on a set of parameters, which enclose positively invariant regions under a vector field of a given two-dimensional system, and also contains an equilibrium point\cite{HaleKocak}. We believe that it can be very useful in the case of practical optimization problems, where choosing the initial condition for the algorithm can be crucial.

We highlight that the ideas to obtain our theoretical results are inspired by the discussions given, on the equation (\ref{original}), in the references \cite{HaleKocak}, Theorem 12.7, and \cite{Ye-Yan}, example 2, pages 8-10.

\noindent
The paper is organized as follows: In section 2 we set the problem and discuss some important aspects related to our setting. Section 3, is devoted to present our main theoretical result, which is to construct a positively invariant region and that is our main theoretical result. We consider pieces of curve which allow to assemble the closed curve that will constitute the boundary of the positively invariant region that is being sought. The section ends with the observation that the instability of the origin allows to conclude the existence of periodic orbits. Section 4 is devoted to recreated numerically our main results. We proposed and algorithm that allows to identify invariant regions and exhibit some picture corresponding to some concrete values of the parameters. Finally, in section 5 we give some concluding remarks where suggest how our results can be useful in determining
the initial condition for the modified vPM's parameters in data assimilation problems, that is, problems associated to the finding the best set of parameters, consistent with a set of data obtained from the observation of the system. There, we emphasize our belief that he strategy presented here can be implemented in other equations.

\section{Setting of the problem, equilibriums and nullclines}

For the purposes of geometrical analysis, as well numerical simulations, we convert the modified van der Pol equation (\ref{lamodificada}) to an equivalent planar system of first-order equations. In fact, by introducing the variables
\begin{equation}
x_1(t)= x(t) ~, ~ x_{2}(t)= \frac{dx}{dt} ~;
\end{equation} 

\noindent
equation (\ref{lamodificada}) becomes
\begin{eqnarray}
\label{ODE}
\dot{x_1} &=&  x_2   \nonumber \\
\dot{x_2} &=& - \alpha ( x_{1}^{2}  - {\nu}^{2} ) x_2 - \frac{x_1 (x_1 + d ) ( x_1 + e )}{ed} ~.
\end{eqnarray}

Our problem is to study the system (\ref{ODE}) and the main goal is to construct a positively invariant region, depending on the parameters, containing the origin. After that, our interest is devoted in determined  the existence of periodic orbits and, as application of the theoretical results,  perform numerical simulations.

\noindent
We assume based on \cite{VDPM} , \cite{Zduniak} and  \cite{Acosta} that
\begin{eqnarray}
\label{Conditions}
\alpha >0 & \mbox{and} & 0 < \nu < e < d \leq 2e. \label{condiciones} 
\end{eqnarray}

The solutions of (\ref{ODE})  that are constant are defined through the equilibrium points:
$$ P_{1}(0 , 0)~, ~ P_{2}(-e , 0) ~, ~ P_{3}(-d , 0)~. $$
The analysis of the stability for those points leads to the following scenarios:
\begin{itemize}
\item
$P_{1}$ is unstable, for each  configuration of parameters $\alpha$, $\nu$, $e$ and $d$ that satisfies (\ref{Conditions}).
\item
$P_{2}$ is unstable, particularly a saddle point, for each configuration of parameters $\alpha$, $\nu$, $e$ and $d$ that satisfies (\ref{Conditions}).
\item
$P_{3}$ is, under the assumption $ \alpha >0$, $ 0< \nu < e <d  $, asymptotically stable. 
\end{itemize}
The nullclines are defined by the curves where either $ \dot{x_1}=0$ or  $ \dot{x_2}=0$ and allow to indicate where the flow is vertical or horizontal. For our system the flow is vertical on the $x_1$-axis and it is horizontal for those points which satisfy the equation 
\begin{eqnarray}~
\label{N}
- \alpha ( x_{1}^{2}  - {\nu}^{2} ) x_2 - \frac{x_1 (x_1 + d ) ( x_1 + e )}{ed} =0~.
\end{eqnarray}
The graph corresponding to (\ref{N})  has three components, two vertical asymptotes with equations $ x_{1} = \pm {\nu} $ and one oblique asymptote with equation
\begin{eqnarray}
\label{AO}
x_1 + \alpha edx_2 + e+ d=  0 ~.
\end{eqnarray}
\hspace{1mm}
\begin{figure}
\begin{center}
	\includegraphics[width=13cm, height=7cm]{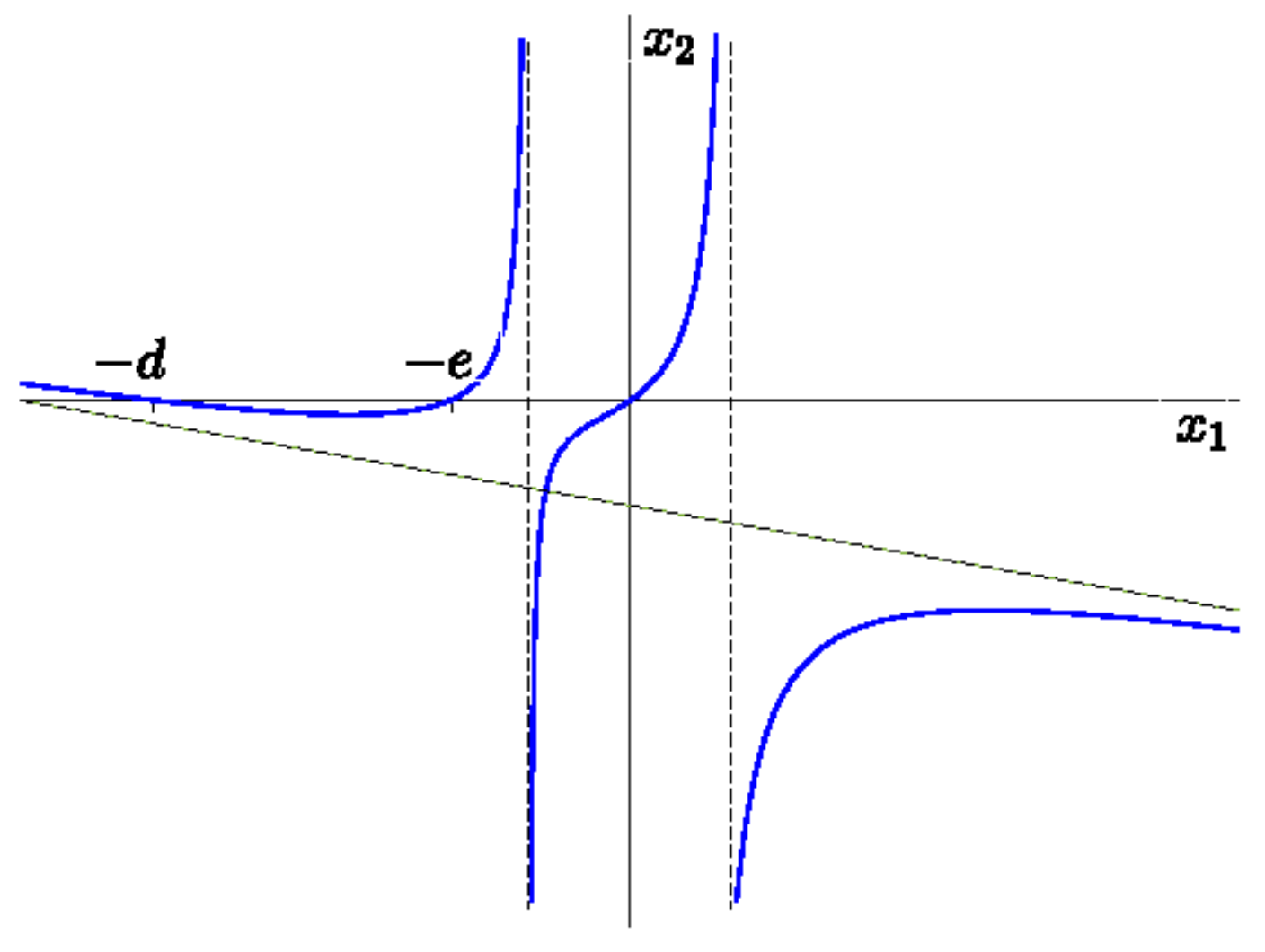}
	\label{Fig:Nullcline}
	\caption{ Nullclines of systems (\ref{ODE}).}
\end{center}
\end{figure}

\section{Existence of positive invariant regions and  periodic orbits}

\begin{figure}
\begin{center}
	\includegraphics[width=13cm, height=7cm]{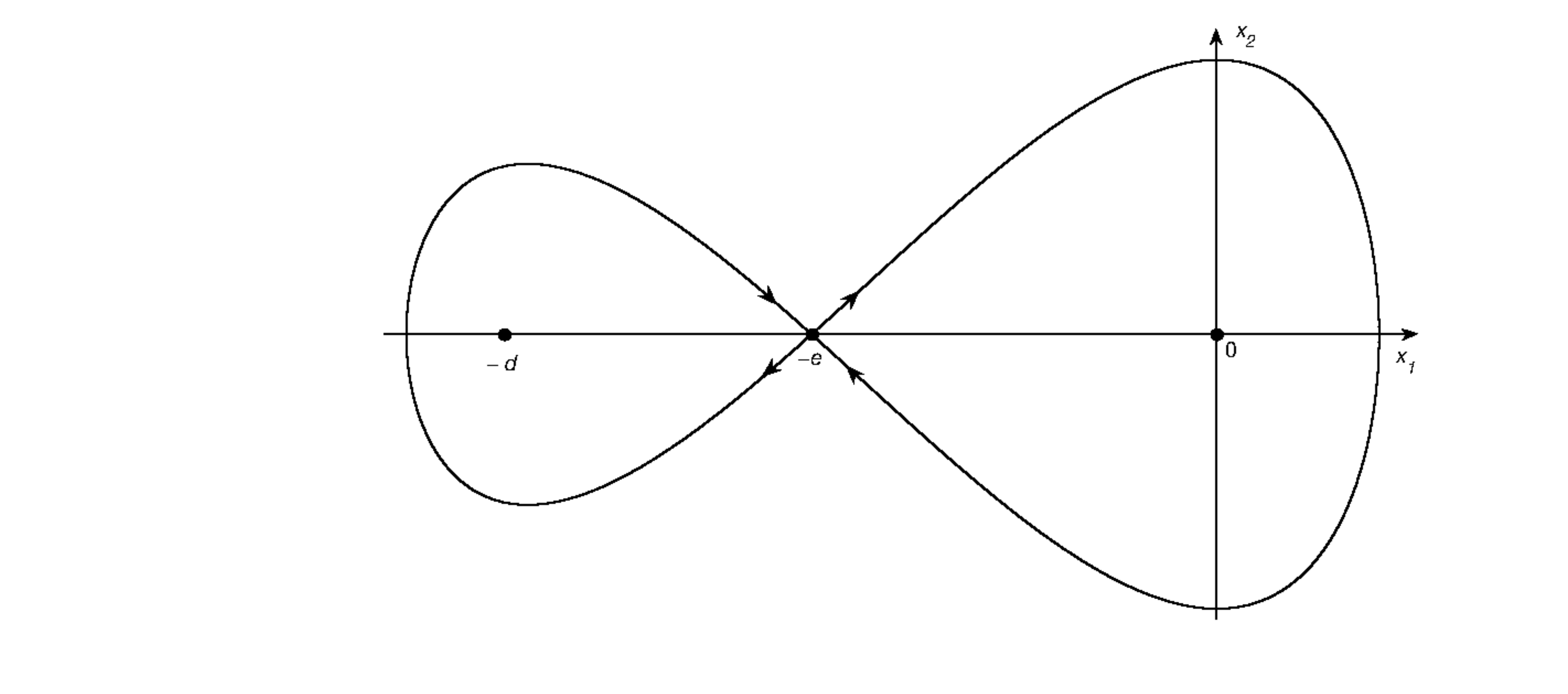}
	\caption{Graph corresponding to (\ref{CQP}).}
	\label{Lazos}
\end{center}
\end{figure}

Our main task, in this section, is to construct a positively invariant region bounded by a closed curve $\cal K$ that will encircles the origin.The curve $\cal K$ will be obtained  by piecing together various curve segments which are coming from some systems that result when a few  terms  of the system (\ref{ODE}) are neglected, the nullclines of (\ref{ODE}) and some  line segments. Specifically, the systems are

\begin{eqnarray}
\dot{x_1} &=&  x_2   \nonumber \\
\dot{x_2} &=& - x_1 + \alpha {\nu}^{2}x_2  \nonumber \\
\dot{x_1} &=&  x_2   \nonumber \\
\dot{x_2} &=& - x_1  \nonumber
\end{eqnarray}

\noindent
and
\begin{eqnarray}
\label{CS}
\dot{x_1} &=&  x_2   \nonumber \\
\dot{x_2} &=&  - \frac{x_1 (x_1 + d ) ( x_1 + e )}{ed} ~.
\end{eqnarray}
Before starting to build the curve $\cal K$, we point out some important aspects of the system (\ref{CS}). First of all,  systems (\ref{ODE}) and (\ref{CS}) share the same equilibrium points, i.e, $ P_{1}(0 , 0)$,$ P_{2}(-e , 0) $ and $ ~ P_{3}(-d , 0) $. Moreover, for (\ref{CS}) it can be shown that $P_2$ is a saddle point, while $P_1$ and $P_3$ are centers.
Also system (\ref{CS}) is conservative  and 
\begin{eqnarray*}
E (x_1 , x_2 ) = \frac{1}{2}x_{2}^{2}+ x_{1}^{2} \left(  \frac{1}{2} + \frac{1}{3}\left( \frac{1}{e} + \frac{1}{d}\right)  x_1 + \frac{1}{4ed} x_{1}^{2} \right) 
\end{eqnarray*}
is a conserved quantity. Thus, given an initial data $(x_{1}^{0}, x_{2}^{0})$ for (\ref{CS}), we have that the corresponding orbit lies on the graph associated to the equation 
\begin{eqnarray}
\label{CQ}
x_{2}^{2} = 2 E\left( x_{1}^{0}, x_{2}^{0} \right) - x_{1}^{2} \left(  1 + \frac{2}{3}\left( \frac{1}{e} + \frac{1}{d}\right)  x_1 + \frac{1}{2ed} x_{1}^{2} \right) ~.
\end{eqnarray}
Particularly we are interested when $x_{1}^{0}=-e$ and  $ x_{2}^{0}=0$, i.e, the saddle point point $P_{2}$. In this scenario the corresponding equation is
\begin{eqnarray}
\label{CQP}
x_{2}^{2} = \frac{1}{3}e^{2}\left( 1 - \frac{e}{2d} \right)  - x_{1}^{2} \left(  1 + \frac{2}{3}\left( \frac{1}{e} + \frac{1}{d}\right)  x_1 + \frac{1}{2ed} x_{1}^{2} \right) ~
\end{eqnarray}
and the graph associated to (\ref{CQP}), see figure \ref{Lazos},  contains two loops, one of them enclosed the point $P_1$ and the other enclosed the point $P_3$. 

Now, we start building the curve $\cal K$. The first piece is contained in the nullclines, corresponding to $\dot{x_{2}}=0$ . It starts at the equilibrium point $P_{2}(-e,0)$ and ends at a certain point $A(x_{10}, x_{20})$, with $ -e < x_{10} < - \nu $. To obtain precisely the coordinates of $A$ we study the points of tangency of an orbit of the system

\begin{eqnarray}
\label{ODE1}
\dot{x_1} &=&  x_2  \nonumber  \\
\dot{x_2} &=& - x_1 + \alpha {\nu}^{2}x_2  
\end{eqnarray}

\noindent
with the component of the set of points $(x_1 , x_2)$ satisfying  the equation (\ref{N}) and such that  $x_1 < - \nu$.
By implicitly differentiating equation (\ref{N}) it is obtained
$$-2 \alpha x_1 x_2 \dot{x_1}- \alpha \left( x_{1}^{2} - {\nu}^{2}\right)\dot{x_2}-\dot{x_1}- \frac{2(e+d)}{ed}x_1 \dot{x_1}- \frac{3}{ed}x_{1}^{2}\dot{x_1}=0 ~,$$
and using (\ref{ODE1}) lead to the equation
$$ -2 \alpha x_1 x_{2}^{2} - \alpha \left( x_{1}^{2} - {\nu}^{2}\right)\left( -x_1 +  \alpha {\nu}^{2}x_2 \right)-x_2 - \frac{2(e+d)}{ed}x_1 x_2 - \frac{3}{ed} x_{1}^{2}x_2 =0 ~.$$

Now, the fact that the points $ (x_1 , x_2)$ satisfy (\ref{N}) allows, after we realize some computations, to obtain the equation 
\begin{equation}
\label{Etangency}
-2 x_{1}^{2}\eta^{2}  + \alpha ed \xi \left( \alpha {\nu}^{2} \xi -1 \right) \eta 
+ \left( 3x_{1}^{2} +2(e+d)x_1  \right) \xi \eta   
+ {\alpha}^{2}e^{2}d^{2} {\xi}^{3}=0 ~.
\end{equation}

\noindent
where, $\eta = \left( x_1 + e \right)\left( x_1 + d \right)$ and $\xi = \left(x_{1}^{2} - {\nu}^{2} \right)$. 

When $x_1= -e$, the left-hand side of the equation (\ref{Etangency}) is positive, and when $x_1$ is close to $- \nu$ the left-hand side is negative; hence, the equation (\ref{Etangency}) has solutions in the interval $ (-e , - \nu )$. Among those solutions, the one nearest to $- \nu$ is the one we take as $x_{10}$, the first coordinate of $A$. The set of points corresponding to the equation (\ref{N}) between $P_2$ and $A$ defined the first piece of $\cal K$, we denote it by $P_2 A$. It is clear that the orbits of the equation (\ref{ODE}) are crossing $P_2 A$ from left to right.

\noindent
For the second piece of $\cal K$, starting at $A$ we follow the orbit of the system (\ref{ODE1}) until it hits the line of the equation
\begin{eqnarray}
\label{recta1}
-x_1 + \alpha {\nu}^{2}x_2 =0 ~.
\end{eqnarray}
It produces a point $B(x_{11}, x_{21})$. The second piece of $ \cal K$, which is denoted by $AB$, is obtained as the piece of the orbit of (\ref{ODE1}) that starts at $A$ and ends at $B$. Since 
$$-x_1 + \alpha {\nu}^{2}x_2 >  - \alpha ( x_{1}^{2}  - {\nu}^{2} ) x_2 - \frac{x_1 (x_1 + d ) ( x_1 + e )}{ed} $$
for the points of $AB$, except possible when $x_1=0$, the orbits of the equation (\ref{ODE}) are crossing the curve $AB$ from left to right.

\noindent
For the third piece of our curve, we consider two  cases: $x_{11} \geq \nu$ and $ x_{11}< \nu$. At this moment we pay attention to the first case, the other will be discuss much later. Now, starting at $B$  we follow the orbit of the system
\begin{eqnarray}
\label{ODE2}
\dot{x_1} &=  x_2   \\
\dot{x_2} &= - x_1  
\end{eqnarray}
until it hits the $x_1$-axis. It produces a point $C(x_{12}, x_{22})$, being $x_{22}=0$. The arc of circumference $BC$ is the third piece of $\cal K$. Since
$$  - \alpha ( x_{1}^{2}  - {\nu}^{2} ) x_2 - \frac{x_1 (x_1 + d ) ( x_1 + e )}{ed} < -x_1 $$
for the points of $BC$, the orbits of the equation (\ref{ODE}) are crossing the curve $BC$ from left to right.

\noindent
To obtain the fourth piece of our curve we consider the vertical line of equation $x_1 = x_{12}$. The intersection of this line with the line of equation (\ref{AO}) is a point $D(x_{13} , x_{23})$, being $x_{13}=x_{12}$. The fourth piece of $\cal K$ is $CD$ and clearly the orbits of the equation (\ref{ODE}) are crossing this segment line from right to left.

\noindent
For the fifth piece of $ \cal K$, starting at $D$ we follow the orbit of (\ref{ODE1}) until it hits the line of equation (\ref{recta1}) and determines a point $E(x_{14}, x_{24})$. The curve $DE$ is the fifth piece of $ \cal K$. Now, since
$$  - \alpha ( x_{1}^{2}  - {\nu}^{2} ) x_2 - \frac{x_1 (x_1 + d ) ( x_1 + e )}{ed} > -x_1 + \alpha {\nu}^{2}x_2 $$
for the points of $DE$, the orbits of the equation (\ref{ODE}) are crossing the curve $DE$ from right to left.
$$ {}$$
Before continuing building the remaining pieces of $ \cal K$ we  now select the parameters $ \alpha , \nu , e $ and $d$ in such a way that the point $E$ is inside the region enclosed by the set of points that satisfies the equation (\ref{CQP}) with the condition $-e < x_{1} \leq - \nu$. Thus, we impose the following two conditions on the coordinates of point $ E$
\begin{eqnarray}\label{Eq:Condition_1}
-e < x_{14} \leq -\nu
\end{eqnarray}
and
\begin{eqnarray}\label{Eq:Condition_2}
- \sqrt{\frac{1}{3}e^{2}\left( 1 - \frac{e}{2d} \right)  - x_{14}^{2} \left(  1 + \frac{2}{3}\left( \frac{1}{e} + \frac{1}{d}\right)  x_{14} + \frac{1}{2ed} x_{14}^{2} \right)} \leq x_{24} <0~.
\end{eqnarray}
$$ {}$$
The sixth piece of $\cal K$ is the vertical line segment going from the point $E$ and ending in the point $F(x_{15}, x_{25})$, being $x_{15}=x_{14}$, which is obtained when we set in the equation  (\ref{CQP}) $x_1= x_{14}$. Thus

$$
x_{25} =- \sqrt{\frac{1}{3}e^{2}\left( 1 - \frac{e}{2d} \right)  - x_{14}^{2} \left(  1 + \frac{2}{3}\left( \frac{1}{e} + \frac{1}{d}\right)  x_{14} + \frac{1}{2ed} x_{14}^{2} \right)} ~.
$$
The orbits of the equation (\ref{ODE}) are crossing the vertical line segment $EF$ from right to left.

\noindent
Now, under the assumption $ x_{11} \geq \nu $, we complete  the last piece of $ \cal K$. In fact, starting  at the point $F$ we follow the points that satisfy equation (\ref{CQP}) until the equilibrium point $P_{2}(-e,0)$ is reached and the arc of curve $FP_2$ is the last piece of $ \cal K$.  Since
$$  - \alpha ( x_{1}^{2}  - {\nu}^{2} ) x_2 - \frac{x_1 (x_1 + d ) ( x_1 + e )}{ed} >  - \frac{x_1 (x_1 + d ) ( x_1 + e )}{ed} $$
for the points of $FP_2$, the orbits of the equation (\ref{ODE}) are crossing the curve $FP_2$ from right to left.

\noindent
Finally all the previous discussion allows us to conclude that the region encircled by the closed curve ${\cal K} = P_2ABCDEFP_2$ is positively invariant for the flow of the system (\ref{ODE}).

\noindent
We conclude the discussion about the existence of a positively invariant region considering the assumption $x_{11} < \nu$. In this scenario our closed curve will have eight pieces. We insert the horizontal line segment with extreme points $B(x_{11}, x_{21})$, $B_{1}( \nu , x_{21})$ and observe that the orbits of the equation (\ref{ODE}) are crossing the segment from top to bottom. 
The rest of the discussion continues as before and so our curve is now ${\cal K} = P_2ABB_{1}CDEFP_2$.

\noindent
Concerning to the existence of periodic orbits we can now affirm, due that the origin is an unstable focus, that the region enclosed by $ \cal K$ contains periodic orbits.
\bigskip

\section{Numerical approach to the positively invariant region }
\begin{figure}
\centering
\includegraphics[width=10cm, height=7cm]{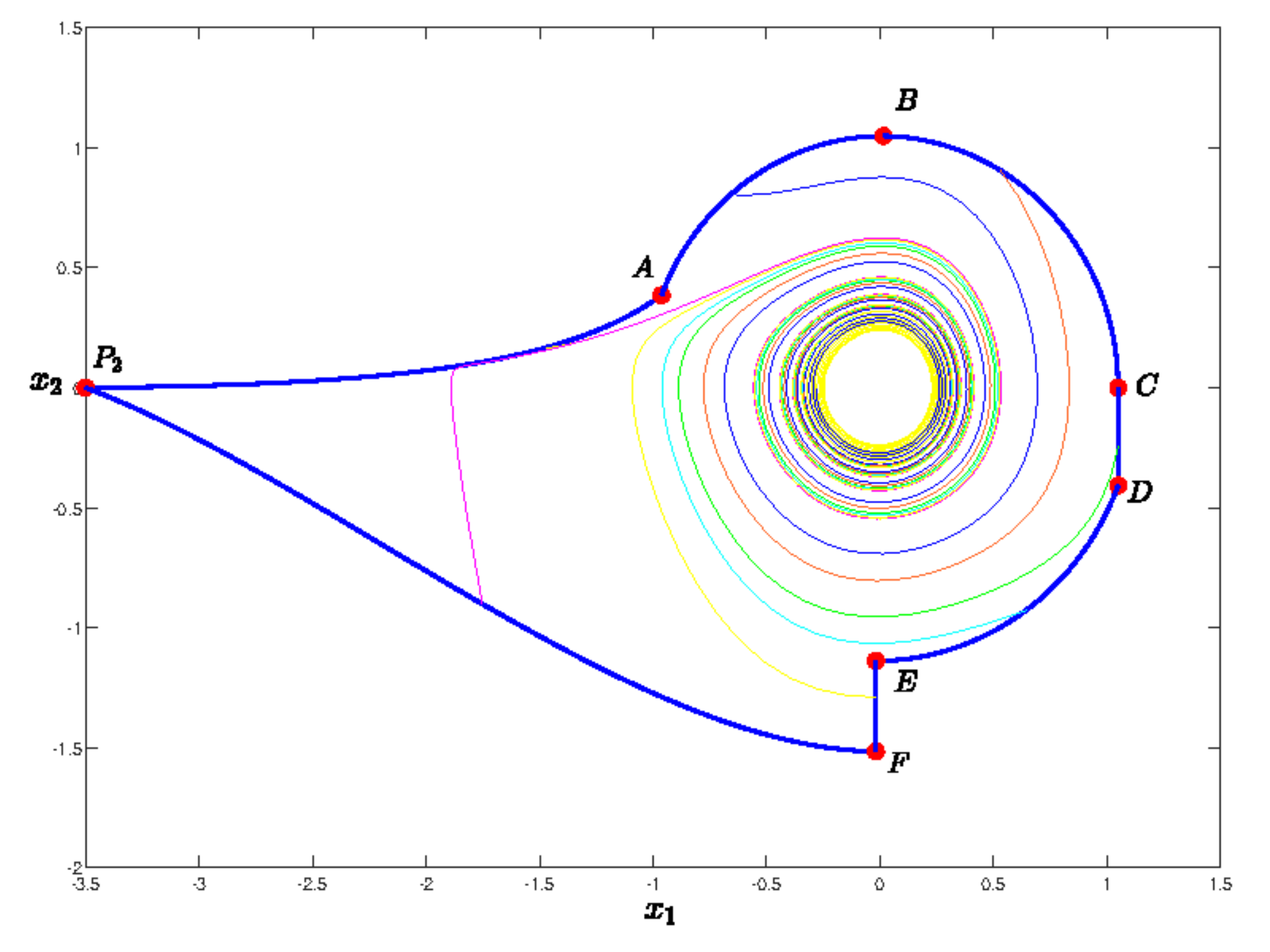}
\caption{The thick line shows the estimate of the invariant region and six orbits generated from initial conditions on the border of the invariant region.}
\label{fig2:Invariant_Region_And_Solutions.eps}
\end{figure}

In this section, we present some numerical results to support our theoretical analysis. 
Specifically, an algorithm to construct the curve $\cal K$  that encloses the invariant region for certain values of the parameters $\alpha$, $e$, $\nu$ and  $d$. The numerical results show that orbits of the system (\ref{ODE})  starting at the $ \cal K$-curve converge to the periodic orbit,
Figure \ref{fig2:Invariant_Region_And_Solutions.eps}. Also, orbits that start outside the positively invariant region are shown. In this case it is observed that they can enter the positively invariant region or move  away,  Figure \ref{fig:Curva_Encierra_Region_Invariante_Orbitas_Externas_2}.
Moreover, although we do not discuss that in this work, our numerical simulations indicate that for each configuration of the parameters with which a positively invariant region is obtained, it contains a single periodic orbit that encloses the origin.

Next, an algorithm is proposed to construct the curve that encloses the positively invariant region obtained in the previous section, for a specific set of parameters $\alpha$, $e$, $\nu$, $d$:

\begin{figure}
\centering
\includegraphics[scale=0.5]{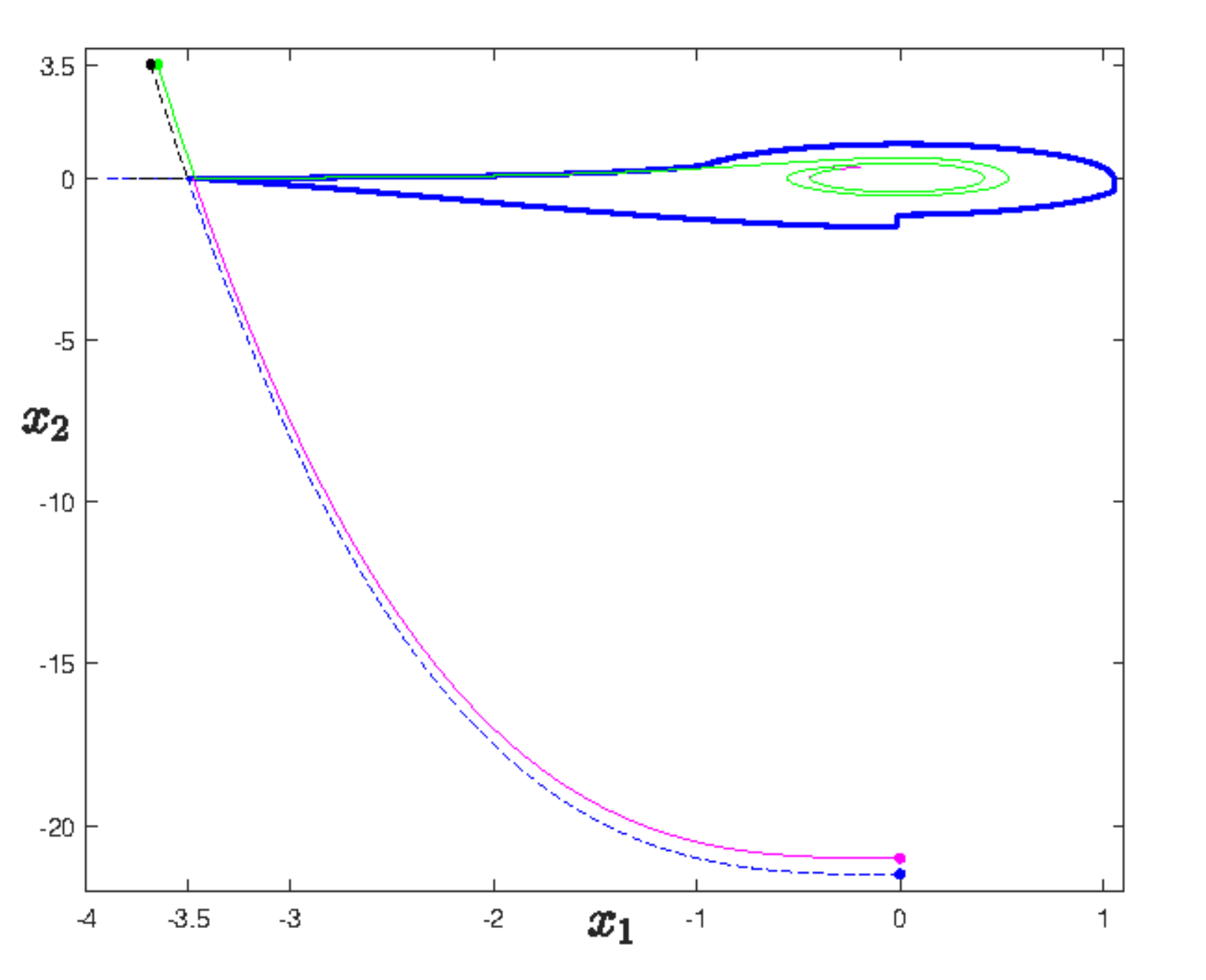}
\caption{The thick line shows the estimate of the invariant region presented in the previous figure and four orbits generated from initial conditions selected out of the invariant region.}
\label{fig:Curva_Encierra_Region_Invariante_Orbitas_Externas_2}
\end{figure}

\begin{enumerate}
\item To obtain the tangency point $A(x_{10}, x_{20})$, of the nullcline and the one orbit of the system (\ref{ODE1}) by solving  (\ref{Etangency}) in  the  interval  $ (-e , - \nu )$, and choosing the nearest solution to $- \nu$. 
\item To obtain $P_2A$,  the first piece  of  $\cal K$. This  is the set of points   corresponding to the equation (\ref{N}) between the points $P_2(-e,0)$ and $A$.
\item To obtain $AB$, the second piece of $\cal K$. It is  the orbit of the system  (\ref{ODE1}) starting at $A$ until  it intersects the line of equation (\ref{recta1}) at $B(x_{11},x_{21})$.
\item If $x_{11} \geq \nu$, obtain the arc of circumference $BC$. It is the third piece of $\cal K$. It is a piece of  the orbit of the system (\ref{ODE2}) starting at $B$ until at point $C(x_{12},0)$. Else, do the steps  4.1-4.2
\begin{enumerate}
	\item[4.1] Inserting the horizontal line segment with extreme points $B(x_{11}, x_{21})$, $B_{1}( \nu , x_{21})$
	\item [4.2] Obtain the arc of circumference $B_{1}C$.  It is a piece of  the orbit of the system (\ref{ODE2}) starting at $B_1$ until at point $C(x_{12},0)$. 
\end{enumerate}

\item To Obtain the piece $CD$, the next piece of $\cal K$. It is the vertical line $x_1=x_{12}$, which starting at $C$ and  ends at the point $D(x_{13},x_{23})$. The point $D$ is the intersection of this line with the line (\ref{AO})
\item Next, find the curve $DE$. This curve starts at $D$ and  follows the orbit of (\ref{ODE1}) until it hits the line of equation (\ref{recta1})  at the point $E(x_{14}, x_{24})$. Here, the values of the parameters $\alpha, \nu, e$ and $d$, must be chosen in such a way that  $x_{14}$ satisfies condition (\ref{Eq:Condition_1}) and $x_{24}$ satisfies condition (\ref{Eq:Condition_2}).

\item Find the curve $FP_2$, with $F(x_{15}, x_{25})$.  This curve is the set of points that satisfy equation (\ref{CQP})  with  $x_2\leq 0$ and $-e\leq x_1\leq x_{14}$.  The first coordinate of  $F$ is  $x_{15}=x_{14}$, an the second coordinate of the point  $F$   is obtained by replacing $x_1= x_{14}$ in the equation  (\ref{CQP}).

\item The next piece of $ \cal K$, is  the vertical line segment $EF$. This piece completes the closed curve  $ \cal K$.
\end{enumerate}

In the example presented in this article, the parameters where chosen as $\alpha=1.5$,~ $\nu=0.1$,~ $e=3.5$,~ $d=4.0$ and the numerical implementation was carried out in Mathlab, using the library ODE45, with parameters: error tolerance$=10^{-13}$, absolute error tolerance $= 10^{-12}$ and norm control by default.

\bigskip
\section{Conclusion}
As final comments, we will offer some ideas that will allow us to frame the strategy presented, in a general classification scheme, for the  particular case  of electrocardiographic (ECG) signals. Although we believe that it can be used in more general cases.

In the case that concerns us, the classification problem, colloquially speaking, consists in deciding whether an ECG signal corresponds to a healthy patient or not. The computational cost of this classifier design depends strongly on the size of the space where the electrocardiograms live, that is, the number of data in each one. One alternative to reduce this dimension is matching the ECG signal to a set of four vdPM parameters.

Our approach, associates an invariant region (the region containing the orbit) to a $4$-tuple representing the system. In this way, we contribute in the direction of finding a relationship between parameters and periodic orbits, that can help to reducing the computational cost of the classification algorithm.

The design of an algorithm that makes it possible to find, from our scheme, an optimal set of parameters that generates the invariant region that contains a particular ECG, is an open problem and the object of our future research.
\bigskip

\section{Acknowledgment} 
\noindent
A. Acosta thanks to Instituci\'on Universitaria Pascual Bravo in Medell\'in, Colombia, for its support in this work. 
P. García thanks to Dr. O. Rojas and the members of the Department of Computer Applications in Science and Engineering of the Barcelona Supercomputer Center (BSC), for their support during part of the development of this work.

\bibliographystyle{plain} 
\bibliography{references}

\end{document}